\documentclass{tac}

\usepackage{amsmath}
\usepackage{mathtools}
\usepackage{latexsym,amssymb}
\usepackage{stmaryrd}
\usepackage{graphicx}
\usepackage{tikz-cd}
\usepackage{mathpartir}
\usepackage[mathscr]{euscript}

\usetikzlibrary{calc,decorations.pathmorphing}
\tikzset{curve/.style={settings={#1},to path={(\tikztostart)
    .. controls ($(\tikztostart)!\pv{pos}!(\tikztotarget)!\pv{height}!270:(\tikztotarget)$)
    and ($(\tikztostart)!1-\pv{pos}!(\tikztotarget)!\pv{height}!270:(\tikztotarget)$)
    .. (\tikztotarget)\tikztonodes}},
    settings/.code={\tikzset{quiver/.cd,#1}
        \def\pv##1{\pgfkeysvalueof{/tikz/quiver/##1}}},
    quiver/.cd,pos/.initial=0.35,height/.initial=0}
\tikzset{tail reversed/.code={\pgfsetarrowsstart{tikzcd to}}}
\tikzset{2tail/.code={\pgfsetarrowsstart{Implies[reversed]}}}
\tikzset{2tail reversed/.code={\pgfsetarrowsstart{Implies}}}
\tikzset{no body/.style={/tikz/dash pattern=on 0 off 1mm}}

\newcommand*{\UU}{\mathcal{U}}
\newcommand*{\E}{{\mathscr{E}}}
\newcommand*{\sE}{\mathrm{s}\mathscr{E}}
\newcommand*{\inftyone}{\ensuremath{(\infty,1)}}
\DeclareMathOperator\Op{op}
\newcommand*{\Simplex}{{\text{\raisebox{-0.1ex}{\includegraphics[height=2ex]{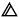}}}}}
\newcommand{\totalty}[1]{\widetilde{#1}}
\newcommand*{\fibarr}{\twoheadrightarrow}
\newcommand*{\cofibarr}{\rightarrowtail}
\newcommand{\ExtAt}[1]{{\mathrm{Ext@}#1}}
\newcommand{\lambdaAt}[1]{{\lambda\mathrm{@}#1}}
\newcommand{\LambdaAt}[1]{{\lambdaAt{#1}}}
\newcommand{\appAt}[1]{{\mathrm{app@}#1}}
\newcommand{\app}{{\mathrm{app}}}
\newcommand{\exten}[4]{\left\langle\mathchoice{\textstyle\prod_{#1}}{\textstyle\prod_{#1}}{\scriptstyle\prod_{#1}}{\scriptscriptstyle\prod_{#1}} #2 \middle|^{#3}_{#4}\right\rangle}
\newcommand{\ndexten}[4]{\left\langle #1 \to #2 \middle|^{#3}_{#4}\right\rangle}
\newcommand{\ccexten}[4]{\langle\Pi[#1][#2] |^{#3}_{#4}\rangle}
\newcommand{\sh}[2]{\{#1\mid #2\}}
\newcommand{\shape}{\;\mathrm{shape}}
\newcommand{\ctx}{\;\mathrm{ctx}}
\newcommand{\types}{\vdash}
\newcommand{\type}{\;\mathrm{type}}
\newcommand{\lam}[1]{\lambda #1.\,}
\newcommand*{\ie}{\emph{i.e.}}
\newcommand*{\cf}{cf.}
\newcommand*{\Cf}{Cf.}
\newcommand*{\eg}{e.g.}
\newcommand*{\unit}{\mathbf{1}}
\DeclareMathOperator\id{id}
\DeclarePairedDelimiter\sem{\llbracket}{\rrbracket}
\DeclarePairedDelimiter\angled{\langle}{\rangle}
\DeclarePairedDelimiterX\pair[2]\langle\rangle{#1,#2}
\newtheoremrm{rem}{Remark}
\newtheoremrm{notation}{Notation}
\newtheoremrm{constr}{Construction}

\newcommand{\jdeq}{=}
\renewcommand{\hookrightarrow}{\cofibarr}

\usepackage[colorlinks=true]{hyperref}
\hypersetup{allcolors=[rgb]{0.1,0.1,0.4}}
\author{Jonathan Weinberger}

\thanks{I am grateful to the Fowler School of Engineering at Chapman University for generous support of this work. I am particularly thankful to the Fletcher Jones Foundation and
their award of a Fletcher Jones Foundation Faculty Fellowship in
Engineering '25--'28 and the ensuing generous funding of parts of this work.
I also thank the Schmid School of Science and Technology as well
as the Center of Excellence in Computation, Algebra, and Topology
(CECAT), both at Chapman University, for my affiliation with them. For financial support during some stages of this project, I gratefully acknowledge the support of the US Army Research Office under MURI Grant W911NF-20-1-0082 and of the Max~Planck-Institute~for~Mathematics~(MPIM), Bonn, Germany. MPIM has also provided great hospitality. I am also grateful to the Hausdorff Research Institute for Mathematics in Bonn, Germany, for hosting us as part of the trimester ``Prospects of formal mathematics,'' funded by the Deutsche Forschungsgemeinschaft (DFG, German Research Foundation) under Germany's Excellence Strategy – EXC-2047/1 – 390685813. I thank Thomas Streicher and Ulrik Buchholtz for useful discussions and explanations. For valuable feedback and discussions I am grateful to Mathieu Anel, Steve Awodey, Sina Hazratpour, Tim Campion, Nima Rasekh, Emily Riehl, Maru Sarazola, Michael Shulman, and Dominic Verity. Nima Rasekh also pointed out and corrected some references. Furthermore, I am grateful to the anonymous referee, Emily Riehl, and Mike Shulman for many helpful comments to improve the paper.}

\address{
Fowler School of Engineering,
Schmid College of Science and Technology,
 \\
Chapman University, \\[5pt]
Swenson Hall, 1 University Dr, Orange, CA 92866, USA}

\title{Strict stability of extension types}

\copyrightyear{2026}


\eaddress{jweinberger@chapman.edu}

\begin{document}

\maketitle
\begin{abstract}
	The theory of $\inftyone$-categories can be developed synthetically in an augmentation of homotopy type theory introduced by Riehl--Shulman. Central to their development is an additional type forming operation called extensions. The original article sketches the semantics of this formal system, explaining how the simplicial homotopy theory can be used to reason about $\inftyone$-categories presented using the Segal space model. However, they leave it open to demonstrate the strict stability of extension types. We prove this using the splitting method of Voevodsky, later generalized by Lumsdaine--Warren to local universes. The practical upshot is that this system has semantics in simplicial objects of an $\infty$-topos, and thus can be used to prove theorems about internal $\infty$-categories in the sense of Martini--Wolf.
\end{abstract}

\section{Introduction}\label{sec:intro}

\subsection{Overview}

In this note, we prove that the \emph{extension types} from Riehl--Shulman's type theory for synthetic $\inftyone$-categories~\cite{RS17}, as explained and further developed in \cite{riehl:2023,riehl:2025,B19,BW21,Wei-intlsums,Wei-2sided,weinberger:chevalley:2024,kudasov:2024} as well as in \cite{ttt,gratzer:2025,gratzer:2026}\footnote{though the latter ones do not make explicit use of strict extension types} can be interpreted in a way so that they are strictly stable under substitution. Semantically, this corresponds to a uniform choice of pullbacks, up to equality of objects and morphisms.

To achieve this, we use the splitting method described by Voevodsky in \cite[Section~4]{VVTySys} and \cite{VVPi}, a self-contained account of which is given by Streicher \cite[Appendix~C]{StrRep,StrFib}. The central technical parts of our treatise as well as the pedagogical recap on the analogous, earlier results for $\Pi$-types mimick the methods used by Awodey in his work on natural models~\cite{Awo18}.

Going back to unpublished work by Lumsdaine--Shulman, extension types play a central role in Riehl--Shulman's \emph{simplicial homotopy type theory}~\cite[Section~2.2]{RS17} where they are used to define hom-types as well as generalizations of those to other shapes. The very important \emph{path types} in cubical type theory~\cite{CCHM2018,OP16} can also be regarded as an instance.

Further applications of extension types are developed in \cite{gratzer2022controlling,zhang2023three}.

\subsection{Contribution}

Our contribution consists in applying a splitting mentioned as in~\cite{voevodsky2014c,LW15,Awo18,StrRep} to the extension type formers in simplicial HoTT so as to make them strictly stable under pullback. This is Theorem~\ref{thm:coh-ext}. We are working with respect to to just one fixed universe. Furthermore, the shapes are taken to be fibrant as well. This simplifies the setting on a technical level while retaining its applicability to the class of models of the form $\sE := \E^{\Simplex^{\Op}}$ for a type-theoretic model topos (TTMT) $\E$, and beyond \cite{rasekh2025}.

An earlier version of the present text is to be found as~\cite[Chapter~6]{jw-phd} of the author's doctoral dissertation under the supervision of Thomas Streicher at TU Darmstadt, Germany.

\subsection{Notation}

Since we work in the extensional type theory of the intended model, we will not distinguish between extensional and intensional equality, and mostly use the symbol $=$ to denote either.

We will often consider type families $A: \Gamma \times \psi \to \UU$ and their reindexings either along the shape part of the context via a map $j : \varphi \hookrightarrow \psi$, or the standard part of the context via a map $\sigma : \Delta \to \Gamma$. To ease the notation, we will supress the part that is not reindexed in the notation and denote both results in the same way, namely as $j^*A := A \circ (\Gamma \times j) : \Gamma \times \varphi \to \UU$ and $\sigma^*A := A \circ (\sigma \times \psi) : \Delta \times \psi$. Analogous conventions will be used for dependent terms.

Sometimes, when considering maps in a comma category 
\[\begin{tikzcd}
	A && B \\
	& \Gamma
	\arrow["\varphi", from=1-1, to=1-3]
	\arrow["\alpha"', from=1-1, to=2-2]
	\arrow["\beta", from=1-3, to=2-2]
\end{tikzcd}\]
we denote them as $\varphi : A \to_\Gamma B$.

\section{Strict stability of $\Pi$-types}

We first recall how to model $\Pi$-types strictly stably under substitution. The appropriate structure will be defined using \emph{universes categories} \`{a} la Voevodsky.

\subsection{Universe categories}

\begin{definition}[Universe categories, \protect{\cite{voevodsky2014c}}]
	Let $\mathscr C$ be a category. A \emph{universe} in $\mathscr C$ is an object $\UU \in \mathscr C$ together with a map $\pi: \widetilde{\UU} \to \UU$ and a choice of pullback squares, assigning to any object $\Gamma \in \mathscr C$ and map $A: \Gamma \to \UU$ maps $p_A: \Gamma.A \to \Gamma$ and $q_A: \Gamma.A \to \widetilde{\UU}$ such that the square
	\[\begin{tikzcd}
		{\Gamma.A} && {\widetilde{\UU}} \\
		\Gamma && \UU
		\arrow["{q_A}", from=1-1, to=1-3]
		\arrow["{p_A}"', from=1-1, to=2-1]
		\arrow["A"', from=2-1, to=2-3]
		\arrow["\pi", from=1-3, to=2-3]
		\arrow["\lrcorner"{anchor=center, pos=0.125}, draw=none, from=1-1, to=2-3]
	\end{tikzcd}\]
	commutes and is cartesian, \ie, a pullback. Any such cartesian square is called a \emph{split cartesian} or \emph{chosen} square.

	A category together with a universe is called a \emph{universe category}.

	We denote the classifying map $\pi$ by a double-headed arrow $\pi : \widetilde{\UU} \twoheadrightarrow \UU$, as well as any pullback of $\pi$.
\end{definition}

We follow Streicher~\cite[Appendix~C]{StrRep,StrFib} and Voevodsky~\cite{VVTySys,VVPi} to define a class of maps strictly stable under substitution. Namely, using meta-theoretic global choice, for each family $A : \Gamma \to \UU$, we select a square:
\[\begin{tikzcd}
	{\Gamma.A} && {\widetilde{\UU}} \\
	\Gamma && \UU
	\arrow["{p_A}"', two heads, from=1-1, to=2-1]
	\arrow["A"', from=2-1, to=2-3]
	\arrow["{q_A}", from=1-1, to=1-3]
	\arrow["\pi", two heads, from=1-3, to=2-3]
	\arrow["\lrcorner"{anchor=center, pos=0.125}, draw=none, from=1-1, to=2-3]
\end{tikzcd}\]
Consider a square:
\begin{equation}
	\begin{tikzcd}
		{\Delta.B} && {\Gamma.A} \\
		\Delta && \Gamma
		\arrow["{p_B}"', two heads, from=1-1, to=2-1]
		\arrow["\sigma"', from=2-1, to=2-3]
		\arrow["{q}", from=1-1, to=1-3]
		\arrow["{p_A}", two heads, from=1-3, to=2-3]
	\end{tikzcd}\label{eq:mor-sq}
\end{equation}
Then, Square~\eqref{eq:mor-sq} is split cartesian if and only if $B = A \circ \sigma$ and $q_B = q_A \circ q$, meaning the squares assemble as follows:
\[\begin{tikzcd}
	{\Delta.B} & {} & {\Gamma.A} && {\widetilde{\UU}} \\
	\Delta && \Gamma && \UU
	\arrow["{p_B}", from=1-1, to=2-1, two heads]
	\arrow["\sigma", from=2-1, to=2-3]
	\arrow["A", from=2-3, to=2-5]
	\arrow["{q}", dashed, from=1-1, to=1-3]
	\arrow["{p_A}", from=1-3, to=2-3, two heads]
	\arrow["{q_A}", from=1-3, to=1-5]
	\arrow["\pi", from=1-5, to=2-5, two heads]
	\arrow["\lrcorner"{anchor=center, pos=0.125}, draw=none, from=1-3, to=2-5]
	\arrow["{q_B}", curve={height=-20pt}, from=1-1, to=1-5]
	\arrow["B"', curve={height=20pt}, from=2-1, to=2-5]
	\arrow["\lrcorner"{anchor=center, pos=0.125}, shift left=5, draw=none, from=1-2, to=2-5]
\end{tikzcd}\]
This then necessarily implies that the left-hand side is a pullback as well.
Since the identites involved are strict equalities between objects and morphisms involving distinguished squares, this models substitution up to equality on the nose, where $\sigma^*A := B$. We write the gap map $q$ as $q = q_{A,\sigma}$.

Both for the definition of $\Pi$- and extension types we will need to consider sections of the maps $p_A : \Gamma.A \twoheadrightarrow \Gamma$. These can be expressed in terms of the universe as follows.

\begin{proposition}[Classification of sections]
	Let $A : \Gamma \to \UU$ be a family and $\Gamma \vdash a : A$ be a section of $A$. Then, in the model there is a correspondence between sections $a : \Gamma \to \Gamma.A$ of $p_A$, \ie,
	\[\begin{tikzcd}
	& {\Gamma.A} \\
	\Gamma & \Gamma
	\arrow["{p_A}", two heads, from=1-2, to=2-2]
	\arrow["a", from=2-1, to=1-2]
	\arrow["{\id_\Gamma}"', Rightarrow, no head, from=2-1, to=2-2]
	\end{tikzcd}\]
	and maps 
	$\overline{a} : \Gamma \to \widetilde{\UU}$ such that:
	\[\begin{tikzcd}
	& {\widetilde{\UU}} \\
	\Gamma & {\UU}
	\arrow["\pi", two heads, from=1-2, to=2-2]
	\arrow["{\overline{a}}", from=2-1, to=1-2]
	\arrow["A"', from=2-1, to=2-2]
	\end{tikzcd}\]
\end{proposition}

\begin{proof}
	Given a section $a : \Gamma \to \Gamma.A$ of $p_A$, we define $\overline{a}$ to be the composite
	\[ \overline{a} : \Gamma \stackrel{a}{\to} \Gamma.A \stackrel{q_A}{\to} \widetilde{\UU}. \]
	On the other hand, given a map $a : \Gamma \to \widetilde{\UU}$ such that $\pi \circ a = A$, we obtain a gap map $\underline{a} : \Gamma \to \Gamma.A$ as in the following diagram:
	\[\begin{tikzcd}
	\Gamma \\
	& {\Gamma.A} & {\widetilde{\UU}} \\
	& \Gamma & \UU
	\arrow["{\underline{a}}"{description}, dashed, from=1-1, to=2-2]
	\arrow["a", curve={height=-12pt}, from=1-1, to=2-3]
	\arrow[curve={height=12pt}, Rightarrow, no head, from=1-1, to=3-2]
	\arrow["{q_A}", from=2-2, to=2-3]
	\arrow["{p_A}"', from=2-2, to=3-2, two heads]
	\arrow["\lrcorner"{anchor=center, pos=0.125}, draw=none, from=2-2, to=3-3]
	\arrow["\pi", two heads, from=2-3, to=3-3]
	\arrow["A"', from=3-2, to=3-3]
	\end{tikzcd}\]
	In particular, as desired $p_A \circ \underline{a} = \id_\Gamma$. By universality of the pullback, one checks that the assignments $\overline{(\cdot)}$ and $\underline{(\cdot)}$ are mutually inverse.
\end{proof}

Next, we make precise how reindexing by a context morphism, both on the level of dependent types and terms corresponds to precomosition of their classifying maps.

\begin{proposition}[Reindexing is precomposition]\label{prop:reindex-is-precomp}
	Let $A : \Gamma \to \UU$ be a family and $a : \prod_\Gamma A$ be a section. Then, reindexing by a morphism $\sigma : \Delta \to \Gamma$ is given by precomposition as follows:
	\[\begin{tikzcd}
	&& {\widetilde{\UU}} \\
	\Delta & \Gamma & \UU
	\arrow["\pi", two heads, from=1-3, to=2-3]
	\arrow["{\sigma^*a}", curve={height=-18pt}, from=2-1, to=1-3]
	\arrow["\sigma"', from=2-1, to=2-2]
	\arrow["{\sigma^*A}"', curve={height=24pt}, from=2-1, to=2-3]
	\arrow["a", from=2-2, to=1-3]
	\arrow["A"', from=2-2, to=2-3]
\end{tikzcd}\]
\end{proposition}

\begin{proof}
	The reindexing of $A$ by $\sigma$ is given by the pullback:
	\[\begin{tikzcd}
	{\sigma^*(\Gamma.A)} & {\Gamma.A} \\
	\Delta & \Gamma
	\arrow["{q_{A,\sigma}}", from=1-1, to=1-2]
	\arrow["{p_{\sigma^*A}}"', two heads, from=1-1, to=2-1]
	\arrow["\lrcorner"{anchor=center, pos=0.125}, draw=none, from=1-1, to=2-2]
	\arrow["{p_{A}}", two heads, from=1-2, to=2-2]
	\arrow["\sigma"', from=2-1, to=2-2]
	\end{tikzcd}\]
	But this can be extended by the classifying pullback for $A$:
	\[\begin{tikzcd}
	{\sigma^*(\Gamma.A)} & {\Gamma.A} & {\widetilde{\UU}} \\
	\Delta & \Gamma & \UU
	\arrow["{q_{A,\sigma}}", from=1-1, to=1-2]
	\arrow["{p_{\sigma^*A}}"', two heads, from=1-1, to=2-1]
	\arrow["\lrcorner"{anchor=center, pos=0.125}, draw=none, from=1-1, to=2-2]
	\arrow["{q_A}", from=1-2, to=1-3]
	\arrow["{p_{A}}", two heads, from=1-2, to=2-2]
	\arrow["\lrcorner"{anchor=center, pos=0.125}, draw=none, from=1-2, to=2-3]
	\arrow["\pi", two heads, from=1-3, to=2-3]
	\arrow["\sigma"', from=2-1, to=2-2]
	\arrow["A"', from=2-2, to=2-3]
	\end{tikzcd}\]
	By composition, we obtain:
	\[\begin{tikzcd}
	{\sigma^*(\Gamma.A)} & {\Gamma.A} & {\widetilde{\UU}} \\
	\Delta & \Gamma & \UU
	\arrow["{{{q_{A,\sigma}}}}", from=1-1, to=1-2]
	\arrow["{{{q_{\sigma^*A}}}}"{description}, curve={height=-24pt}, from=1-1, to=1-3]
	\arrow["{{{p_{\sigma^*A}}}}"', two heads, from=1-1, to=2-1]
	\arrow["\lrcorner"{anchor=center, pos=0.125}, draw=none, from=1-1, to=2-2]
	\arrow["{{{q_A}}}", from=1-2, to=1-3]
	\arrow["\pi", two heads, from=1-3, to=2-3]
	\arrow["\sigma"', from=2-1, to=2-2]
	\arrow["{{{A\sigma}}}"{description}, curve={height=24pt}, from=2-1, to=2-3]
	\arrow["A"', from=2-2, to=2-3]
	\end{tikzcd}\]
	This shows that $\sigma^*A = A \circ \sigma$, which we also write as $\sigma^*A$. 
\end{proof}

Relative to this universe structure, one defines new type formers using appropriate classifying objects, which are also regarded as \emph{generic contexts} for the corresponding type-theoretic rules.

In parallel to their categorical semantics, one can also understand these generic contexts in terms of the internal extensional type theory (ETT) of the ambient topos as is also done in~\cite{KL18,Awo18,LW15}.

We now assume $\mathscr{C}$ to be, moreover, locally cartesian closed category. Given a map $f:B \to A$, we denote the ensuing adjoint triple byp:
\[\begin{tikzcd}
	{\mathscr C/A} && {\mathscr C/B}
	\arrow[""{name=0, anchor=center, inner sep=0}, "{f^*}"{description}, from=1-1, to=1-3]
	\arrow[""{name=1, anchor=center, inner sep=0}, "{\sum_f}"{description}, curve={height=24pt}, from=1-3, to=1-1]
	\arrow[""{name=2, anchor=center, inner sep=0}, "{\prod_f}"{description}, curve={height=-24pt}, from=1-3, to=1-1]
	\arrow["\dashv"{anchor=center, rotate=-90}, draw=none, from=1, to=0]
	\arrow["\dashv"{anchor=center, rotate=-90}, draw=none, from=0, to=2]
\end{tikzcd}\]
With this notation, and conflating objects with their terminal projections, the \emph{polynomial functor associated to a morphism $f:B \to A$}~\cite{gambino2013polynomial,Awo18,aberle:2025} is defined as
the composite:
\[\begin{tikzcd}
	{\mathscr C} && {\mathscr C/B} && {\mathscr C/A} && {\mathscr C}
	\arrow["{\prod_f}", from=1-3, to=1-5]
	\arrow["{\sum_{A}}", from=1-5, to=1-7]
	\arrow["{B^*=(B \times -)}", from=1-1, to=1-3]
	\arrow["{P_f}"{description}, shift right=1, curve={height=24pt}, from=1-1, to=1-7]
\end{tikzcd}\]
The action on objects is given by
\[ P_f(X) := \sum_A \prod_f B^* (X)= \prod_f(B \times X) \in \mathscr C.  \]
Writing $B_a := f^*a$, this reads in the internal type theory of $\mathscr C$ as
\[ P_f(X) = \sum_{a:A} X^{B_a},\]
which shows how the adjoint triple aligns with the internal type theory of $\mathscr C$.

The action on morphisms is given by
\[ P_f(k:Y \to X) = \Big( \prod_f(k \times B): \prod_f(Y \times B) \to \prod_f(X \times B) \Big). \]

For the universal fibration $\pi: \widetilde{\UU} \to \UU$, the total object denotes the generic context
\[ \widetilde{\UU} = \sum_{A:\UU} A =  \sem{A:\UU, a:A},\]
and $\pi : \widetilde{\UU} = \sum_{A:\UU} A \twoheadrightarrow \UU$ acts as projection $\pi(A,a) = A$.

We first recall how to define the structure of $\Pi$-types relative to this structure \`{a} la~\cite[Proposition~2.4]{Awo18}. This will serve to motivate the ideas behind the case of extension types which are a kind of $\Pi$-type with additional judgmental boundary conditions.

\subsection{Strictly stable $\Pi$-types (recollection)}

The formation rule for dependent function types is given in Figure~\ref{fig:pi-form}.
\begin{figure}[h]
	\centering
	\begin{mathpar}
		\inferrule*[right = ($\Pi$-Form)]{\Gamma \types A \type \\ \Gamma.A \types B  \type}{\Gamma \types \prod_{a:A} B(a)}
	\end{mathpar}
	\caption{$\Pi$-types: formation}
	\label{fig:pi-form}
\end{figure}

The premise of this rule is
\[ \Gamma \types A \type \qquad \Gamma.A \types B \type. \]
To split the $\Pi$-type former, we seek to define an object $\UU^\Pi$ which classifies this data in the following sense: for any object $\Gamma$, morphisms $\Gamma \to \UU^\Pi$ correspond to pairs $\pair{A}{B}$ with $A: \Gamma \to \UU$ and $B: \Gamma.A \to \UU$, and this correspondence is natural in $\Gamma$. 

The introduction rule for the $\Pi$-type is given by Figure~\ref{fig:pi-intro}.

\begin{figure}[h]
	\centering
	\begin{mathpar}
		\inferrule*[right = ($\Pi$-Intro)]{\Gamma \types A \type \\ \Gamma.A \types B  \type \\ \Gamma.A \types b:B}{\Gamma \types \lambda a.b:\prod_{a:A} B(a)}
	\end{mathpar}
	\caption{$\Pi$-types: introduction}
	\label{fig:pi-intro}
\end{figure}

To model this rule, we are likewise looking for an object $\UU^\lambda$ such that for any object $\Gamma$, a morphism $\Gamma \to \UU^\lambda$ corresponds to the premises of the introduction rule, \ie, triples $\angled{A,B,b}$ with $A: \Gamma \to \UU$, $B: \Gamma.A \to \UU$, and sections $b:\Gamma.A \to \Gamma.A.B$ of $p_B:\Gamma.A.B \twoheadrightarrow \Gamma.A$. It is well known how to define the desired objects $\UU^\Pi$ and $\UU^\lambda$.

\begin{definition}[Generic contexts for $\Pi$-types]\label{def:gen-ctxt-pi}
Let $\mathscr C$ be a locally cartesian closed category and $\pi: \widetilde{\UU} \to \UU$ a map in $\mathscr C$ (to be thought of as the universe of small types).
\begin{enumerate}
	\item\label{it:gen-ctx-pi-form} The \emph{generic context for $\Pi$-formation} is defined to be the object
	\[ \UU^\Pi := \sum_{A: \UU} \UU^A.\]
	\item\label{it:gen-ctx-pi-intro}  The \emph{generic context for $\Pi$-introduction} is denoted by the object
	\[ \UU^\lambda := \sum_{A: \UU} \widetilde{\UU}^A.\]
\end{enumerate}
\end{definition}

Indeed, one can show that the type $\UU^\Pi$ classifies the data of the premise of the $\Pi$-formation rule, \ie, a dependent type $\Gamma \vdash A$ over some context $\Gamma$, together with a dependent type $\Gamma.A \vdash B$ depending on (the extended context made up from) $\Gamma$ and $A$: 
\[
\Gamma \types A \type \qquad \Gamma.A \types B \type
\]
Analogously, the type $\UU^\lambda$ classifies the data of the premise of the $\Pi$-introduction rule, \ie, dependent types $\Gamma \vdash A$ and $\Gamma.A \vdash B$ as above, together with a dependent term $\Gamma.A \types b:B$ in $B$:
\[
\Gamma \types A \type \qquad \Gamma.A \types B \type \qquad \Gamma.A \types b:B 
\]
This is proven by \cite[Proposition~2.2]{Awo18}:

\begin{proposition}[Generic contexts for $\Pi$-types]
The generic contexts from Definition~\ref{def:gen-ctxt-pi} classify the following data:
\begin{enumerate}
	\item\label{it:pi-form} For any object $\Gamma \in \mathscr C$, maps $\alpha : \Gamma \to \UU^\Pi$ are in natural bijection with pairs $\pair{A}{B}$ where $A: \Gamma \to \UU$ and $B: \Gamma.A \to \UU$.
	\item\label{it:pi-intro}  For any object $\Gamma \in \mathscr C$, maps $\beta: \Gamma \to \UU^\lambda$ are in natural bijection with triples $\angled{A,B,b}$ where $A: \Gamma \to \UU$, $B: \Gamma.A \to \UU$, and $b:\Gamma.A \to \Gamma.A.B$ is a section of $p_B:\Gamma.A.B \to \Gamma.A$.
\end{enumerate}
The postulated bijections are \emph{natural} in the underlying context $\Gamma$, \ie, for any morphism $\sigma: \Delta \to \Gamma$ we have that the composite
\[\begin{tikzcd}
	\Delta && \Gamma && {\UU^\Pi}
	\arrow["\sigma", from=1-1, to=1-3]
	\arrow["\alpha = \pair{A}{a}", from=1-3, to=1-5]
\end{tikzcd}\]
classifies pairs $\pair{A \circ \sigma}{B \circ q_{A,\sigma}}$ arising from the diagram:
\[\begin{tikzcd}
	{\Delta.A\sigma} && {\Gamma.A} && {\widetilde{\UU}} \\
	\\
	\Delta && \Gamma && \UU
	\arrow["{q_{A,\sigma}}"{description}, dashed, from=1-1, to=1-3]
	\arrow["{q_A}"{description}, from=1-3, to=1-5]
	\arrow["{q_{A\sigma}}"{description}, curve={height=-18pt}, from=1-1, to=1-5]
	\arrow["{p_{A\sigma}}"', two heads, from=1-1, to=3-1]
	\arrow["\sigma"{description}, from=3-1, to=3-3]
	\arrow["A"{description}, from=3-3, to=3-5]
	\arrow["\pi", two heads, from=1-5, to=3-5]
	\arrow["{p_A}"', two heads, from=1-3, to=3-3]
	\arrow["A\sigma"{description}, curve={height=18pt}, from=3-1, to=3-5]
	\arrow["\lrcorner"{anchor=center, pos=0.125}, draw=none, from=1-3, to=3-5]
	\arrow["\lrcorner"{anchor=center, pos=0.125}, draw=none, from=1-1, to=3-5]
\end{tikzcd}\]
with $B:\Gamma.A \to \UU$, so that
\[\begin{tikzcd}
	{\Delta.A\sigma} && {\Gamma.A} && \UU
	\arrow["B", from=1-3, to=1-5]
	\arrow["{q_{A,\sigma}}", from=1-1, to=1-3]
	\arrow["{B q_{A,\sigma}}"', curve={height=18pt}, from=1-1, to=1-5]
\end{tikzcd}\]
Given a composite
\[\begin{tikzcd}
	\Delta && \Gamma && {\UU^\Pi}
	\arrow["\sigma", from=1-1, to=1-3]
	\arrow["\beta", from=1-3, to=1-5]
\end{tikzcd}\]
additionally gives a reindexed section arising as follows:
\[\begin{tikzcd}
	{\Delta.A\sigma.a} \\
	& {\Delta.A\sigma.B\sigma} && {\Gamma.A.B} \\
	& {\Delta.A\sigma} && {\Gamma.A}
	\arrow["{\sigma^*a}", dashed, from=1-1, to=2-2]
	\arrow["{a \circ q_{A,\sigma}}", curve={height=-12pt}, from=1-1, to=2-4]
	\arrow[curve={height=12pt}, Rightarrow, no head, from=1-1, to=3-2]
	\arrow["{q_{B,\sigma}}", from=2-2, to=2-4]
	\arrow["{p_{B\sigma}}", two heads, from=2-2, to=3-2]
	\arrow["\lrcorner"{anchor=center, pos=0.125}, draw=none, from=2-2, to=3-4]
	\arrow["{p_B}", two heads, from=2-4, to=3-4]
	\arrow["{q_{A,\sigma}}"', from=3-2, to=3-4]
\end{tikzcd}\]
\end{proposition}

To analyze the generic contexts in the internal language of the ambient lccc, the following well-known result is of central importance.

\begin{proposition}[Distributivity of $\Sigma$ over $\Pi$]\label{prop:dist-law}
	In extensional type theory with judgmental $\Sigma$- and $\Pi$-constructors, let $A$ be a type $A \vdash B$ a dependent type over $A$ and $A.B \vdash C$ a dependent type over the extended context $A.B$. Then
	\[ \prod_{a:A} \sum_{b:B(a)} C(a,b) = \sum_{f : \prod_{x:A} B(x)} \prod_{a:A} C(a,f(a)). \]
\end{proposition}

By Proposition~\ref{prop:dist-law}, we obtain that
\[ \UU^\lambda :=  \sum_{A: \UU } \totalty{\UU}^A = \sum_{A: \UU } \big(A \to \sum_{B:\UU} B\big) = \sum_{A: \UU } \sum_{B:A \to \UU} \prod_{a:A} B(a). \]
This corresponds to the \emph{generic dependent term}, captured by the context
\[ \sem{A : \UU,B : \UU \to A,a:A \vdash b:B(a)}.\]

\begin{figure}[t]
	\centering
	\begin{mathpar}
		\inferrule*[right = ($\Pi$-Elim)]{\Gamma \types A \type \\ \Gamma.A \types B\type  \\ \Gamma \types a : A \\ \Gamma \types f : \prod_{x:A} B(x)}{\Gamma \types \app(f,a):\prod_{a:A} B(a)}
	\end{mathpar}
	\caption{$\Pi$-types: elimination}
	\label{fig:pi-elim}
\end{figure}

\begin{figure}[t]
	\centering
	\begin{mathpar}
		\inferrule*[right = ($\Pi$-$\beta$)]{\Gamma \types A \type \\ \Gamma.A \types B  \type \\ \Gamma \types x : A \\ \Gamma.A \types b:B}{\Gamma \types (\lambda a.b)(x) \jdeq  b[x/a]}
	\end{mathpar}
	\caption{$\Pi$-types: $\beta$-rule}
	\label{fig:pi-beta}
\end{figure}

\begin{figure}[t]
	\centering
	\begin{mathpar}
		\inferrule*[right = ($\Pi$-$\eta$)]{\Gamma \types A \type \\ \Gamma.A \types B  \type \\  \Gamma.A \types b:B}{\Gamma \types b \jdeq \lambda a.b(a)}
	\end{mathpar}
	\caption{$\Pi$-types: $\eta$-rule}
	\label{fig:pi-eta}
\end{figure}

The universe then admits $\Pi$-types (\cf~\emph{$\Pi$-structure}, \cf~\cite[Definition~1.4.2, Theorem~1.4.15]{KL18}), with the standard rules of formation, introduction, elimination, and $\beta$- and $\eta$-computation, if and only if there exist maps $\Pi: \UU^\Pi \to \UU$ and $\lambda: \UU^\lambda \to \widetilde{\UU}$ 
making the square
\[\begin{tikzcd}
	{\UU^\lambda} && {\widetilde{\UU}} \\
	{\UU^\Pi} && {\UU}
	\arrow[from=1-1, to=2-1]
	\arrow["\Pi"', from=2-1, to=2-3]
	\arrow["\lambda", from=1-1, to=1-3]
	\arrow[from=1-3, to=2-3, two heads]
\end{tikzcd}\]
commute, and moreover rendering it as a pullback.

These maps are supposed to implement the formation and introduction rule. So if the universe $\UU$ is indeed closed under $\Pi$-types, we define
\[ \Pi:\UU^\Pi \to \UU, ~ \Pi(A,B) := \prod_A B, \quad  \lambda:\UU^\lambda \to \UU, ~ \lambda(A,B,b) :=\Big \langle \prod_A B, b \Big \rangle,\]
induced from the ambient category being locally cartesian closed. Externally, for maps $p_A: \Gamma.A \twoheadrightarrow \Gamma$, $p_B: \Gamma.A.B \twoheadrightarrow \Gamma.A$ this amounts to defining $\Pi(A,B) := \prod_{p_A}(p_B)$ via the pushforward functor $\prod_{p_A} \vdash (p_A)^*$.\footnote{In a homotopical setting, such as that of \emph{type-theoretic model categories}, the universe $\pi$ classifies \emph{fibrations} $p_A$ which the pushforward is asked to preserve.}

Concretely, we have the following proposition, proven in \cite[Prop.~2.4]{Awo18}.

\begin{proposition}[$\Pi$-type structures]
	\begin{enumerate}
		\item  The formation and introduction rules for $\Pi$-types are modeled by maps
		\begin{align*}
			\Pi & \colon  \UU^\Pi \to \UU, \\
		   \lambda & \colon \UU^\lambda \to \UU. \\
		\end{align*}
		making the following square commute:
		\[\begin{tikzcd}
			{\UU^\lambda} && {\widetilde{\UU}} \\
			{\UU^\Pi} && {\UU}
			\arrow[from=1-1, to=2-1]
			\arrow["\Pi"', from=2-1, to=2-3]
			\arrow["\lambda", from=1-1, to=1-3]
			\arrow[from=1-3, to=2-3]
		\end{tikzcd}\]
		\item The square is a weak pullback, with a distinguished section $\app$ as in
		\[\begin{tikzcd}
			&& {\UU^\lambda} \\
			{\UU^\Pi \times_\UU \widetilde{\UU}} && {\UU^\Pi \times_\UU \widetilde{\UU}}
			\arrow[Rightarrow, no head, from=2-1, to=2-3]
			\arrow[from=1-3, to=2-3]
			\arrow["\app", from=2-1, to=1-3]
		\end{tikzcd}\]
		if and only if the elimination rule with $\beta$-computation holds.
		\item It is a pullback if and only if in addition the $\eta$-computation rule holds.
	\end{enumerate}
\end{proposition}

 In \protect{\cite[Rem.~2.6]{Awo18}} it is shown that modeling the $\Pi$-types via the generic contexts is strictly stable under context substitution:

\begin{proposition}[Strict stability of $\Pi$-types]
	The given interpretation of $\Pi$-types is strictly stable under substitution, \ie, given a context morphism $\sigma : \Delta \to \Gamma$, we have:
		\begin{align*}
			\sigma^*\big(\Pi(A,B)\big) & = \Pi(\sigma^*A, \sigma^*B) \\
			\sigma^*\big(\lambda_{A,B}t.b\big) & = \lambda_{\sigma^*A,\sigma^*B}s.\sigma^*b \\
			\sigma^*\big(f(a)\big) & = (\sigma^*f)(\sigma^*a)
		\end{align*}
\end{proposition}

\section{Strict stability of extension types}

\subsection{General ideas}

The rules for the extension types from~\cite{RS17} are recalled in Appendix~\ref{app:exten}.

We adapt the coherence construction presented in \cite{KL18,LW15,Awo18}, based on \cite{VVTySys,VVPi}, so that in the model extension types \`{a} la \cite{RS17} can be chosen in a way that is strictly stable under pullback.

Our presentation is close to \cite{KL18,Awo18} in style.

In contrast to most of~\cite[Appendix~A]{RS17} we will assume the shapes to be fibrant, which is justified by the intended standard models, namely simplicial objects in an $\infty$-topos.

As in~\cite[Theorem~A.16]{RS17}, in the (injective) model structure of the TTMT presenting $\sE$ we interpret the extension types as ordinary $1$-categorical pullbacks of Leibniz cotensors. Since the model structure at hand is type-theoretic~\cite[Definition~2.12]{ShuInv} the Leibniz cotensors are fibrations, so are all their reindexings, which in particular includes the (dependent) extension types. In the case of non-fibrant shapes, the Leibniz cotensor is still a fibration because the model structure is cartesian monoidal, as discussed in \cite[Section~1, Section~A.2]{RS17}. In particular, the ensuing pullbacks are homotopy pullbacks as proven in \cite[Proposition~A.2.4.4 and Remark~A.2.4.5]{LurHTT}.

In the non-dependent case the extension type $\ndexten{\psi}{A}{\varphi}{a}$ is thus given simply by (ordinary, $1$-categorical) the pullback:
\[\begin{tikzcd}
	{\left\langle \psi \rightarrow A \,\vert^{\varphi}_{a}\right\rangle} && {A^\psi} \\
	{\unit} && {A^\varphi}
	\arrow[from=1-1, to=2-1]
	\arrow["a"', from=2-1, to=2-3]
	\arrow[from=1-1, to=1-3]
	\arrow[from=1-3, to=2-3]
	\arrow["\lrcorner"{anchor=center, pos=0.125}, draw=none, from=1-1, to=2-3]
\end{tikzcd}\]
For the dependent case the construction will have to be relativized as we will see later on.

\begin{notation} We will write $\Pi(A,B)$ and $\ccexten{\varphi}{A}{\varphi}{a}$, respectively, instead of $\prod_A B$ and $\exten{\psi}{A}{\varphi}{a}$, respectively, to designate the (strictly stable) type formers defined as part of the interpretation, in line with~\cite{RS17,KL18,KL18,Awo18}.
\end{notation}

We want to show that the interpretation of the extension types can be chosen in a strictly pullback-stable way. In particular, after~\cite[Definition~A.10, Theorem~A.16]{RS17} given a type family $A \fibarr \Gamma \times \psi$ and a partial section $a:\Gamma \times \varphi \to_{\Gamma \times \psi} A$ substitution of the extension type $\ccexten{\varphi}{A}{\varphi}{a}$ is only considered along (type) context morphisms $\sigma: \Delta \to \Gamma$, leaving the shape inclusion $j: \varphi \cofibarr \psi$ fixed. Hence, we will understand the type former of the extension type as a \emph{family} of type formers $\ExtAt{j}$, given by
\[ \ExtAt{j}(A,a) := \ccexten{\psi}{A}{\varphi}{a}\]
indexed by the shape inclusions $j$. Note that, \emph{a posteriori} this also yields stability with respect to to reindexing along cube context morphisms $\tau:J \to I$ pulling back shape inclusions $j:\varphi \to_I \psi$: By their defining universal property, which only involves a condition on reindexings of the type context $\Gamma$, the (pseudo-stable) extension types from~\cite[Definition~A.10]{RS17} are determined uniquely up to isomorphism. The splitting then yields uniqueness up to equality. Then, considering reindexings along shape maps does not lead out of this class, so the choice remains strictly stable.

We now treat the extension types from \cite{RS17} in a similar fashion.

\subsection{Strictly stable extension types: Generic contexts}
Recall the rules from~\cite[Figure~4]{RS17}.

To form an extension type, we start with a context fibered-over-shapes
\[ \Gamma \fibarr \varphi \cofibarr \Xi\]
and a \emph{separate} shape inclusion
\[ \psi \cofibarr \varphi \cofibarr I.\]

The semantic reason for this is explained right before~\cite[Theorem~A.17]{RS17}. \Cf~also the explanation about the rules at the beginning of~\cite[Section~2.2]{RS17}
	
The extension type is then formed for a pair $\pair{A}{a}$ of a type $A$ and a partial section $a$ where
\begin{equation}\label{eq:input-ext-form}
	(\Gamma \times \psi).A \fibarr \Gamma \times \psi \fibarr \varphi \times \psi \cofibarr \Xi \times I, \quad a: \Gamma \times \varphi \to_{\Gamma \times \psi} A
\end{equation}
In fact, the defining universal property for the extension types \cite[Definition~A.10]{RS17} asks for substitution (pseudo-)stability of $A$ and $a$ along morphisms $\sigma:\Delta \to \Gamma$. Hence, we in fact consider a \emph{family} of type formers $(\ExtAt{j})_{j:\varphi \cofibarr_I \psi}$ indexed externally by the shape inclusions $j:\varphi \to \psi$, so that
\[ (\ExtAt{j})(A,a) = \ccexten{\psi}{A}{\varphi}{a}.\]
Thus, the input data to form the extension type $\ExtAt{j}$ should be represented using a suitable generic context $\UU^{\ExtAt{j}}$ as a morphism
\[ \pair{A}{a}: \Gamma \to \UU^{\ExtAt{j}} \]
with $\pair{A}{a}$ as in~\eqref{eq:input-ext-form}.

Note that we are suppressing the further structure $\Gamma \fibarr \varphi \cofibarr \Xi$ here, in line with~\cite[Definition~A.10]{RS17}.

In analogy to the case of $\Pi$-types, we form the generic contexts as follows. Consider the universe $\pi : \widetilde{\UU} \twoheadrightarrow \UU$ and consider the square induced by exponentiation with the shape inclusion $j: \varphi \hookrightarrow \psi$:
\[\begin{tikzcd}
	{\widetilde{\UU}^\psi} \\
	& {\UU^\psi \times_{\UU^\varphi} \widetilde{\UU}} && {\widetilde{\UU}} \\
	& {\UU^\psi} && {\UU^\varphi}
	\arrow[dashed, from=1-1, to=2-2]
	\arrow["{\widetilde{\UU}^j}", curve={height=-18pt}, from=1-1, to=2-4]
	\arrow["{\pi^\psi}"', curve={height=12pt}, from=1-1, to=3-2]
	\arrow[from=2-2, to=2-4]
	\arrow[from=2-2, to=3-2]
	\arrow["\lrcorner"{anchor=center, pos=0.125}, draw=none, from=2-2, to=3-4]
	\arrow["{\pi^\varphi}", from=2-4, to=3-4]
	\arrow["{\UU^j}"', from=3-2, to=3-4]
\end{tikzcd}\]

To set up our classifcation results we introduce the following notations:
\begin{definition}[Generic contexts for extension types]
	We call
	\begin{enumerate}
		\item the object $\UU^{\ExtAt{j}} := {\UU^\psi \times_{\UU^\varphi} \widetilde{\UU}}$ the \emph{generic context for $\ExtAt{j}$-formation},
		\item the object $\UU^{\LambdaAt{j}} := \widetilde{\UU}^\psi$ the \emph{generic context for $\ExtAt{j}$-introduction},
		\item the object $\UU^{\appAt{j}} := \UU^{\ExtAt{j}} \times_\UU \widetilde{\UU}$ the \emph{generic context for $\ExtAt{j}$-elimination}.
	\end{enumerate}
\end{definition}

These then serve to define our generic contexts, classifying the data needed for the extension type former:

\begin{proposition}[Generic contexts for $\ExtAt{j}$-formation and introduction]\label{prop:ext-class}
The generic contexts classify the following data:
\begin{enumerate}
	\item\label{it:ext-form} $\UU^{\ExtAt{j}}$ classifies families over a shape with partial sections. Concretely: for any object $\Gamma \in \mathscr C$, maps $\alpha : \Gamma \to \UU^{\ExtAt{j}}$ are in natural bijection with pairs $\pair{A}{a}$ where $A: \Gamma \times \psi \to \UU$ and $a: \prod_{\Gamma \times \varphi} A_j$, as defined in the following diagram:
	\[\begin{tikzcd}
	& {\widetilde{\UU}} \\
	{\Gamma \times \varphi} & \UU
	\arrow["\pi", two heads, from=1-2, to=2-2]
	\arrow["a", from=2-1, to=1-2]
	\arrow["{j^*A}"', from=2-1, to=2-2]
	\end{tikzcd}\]
	\item\label{it:ext-intro} $\UU^{\ExtAt{j}}$ classifies families over a shape with (total) sections, which in turn correspond to partial sections together with an extension and witnessing equality. Concretely: for any object $\Gamma \in \mathscr C$, maps $\beta: \Gamma \to \UU^{\lambdaAt{j}}$ are in natural bijection with triples $\angled{A,a,b,\_}$ where $A: \Gamma \times \psi \to \UU$, $a:\prod_{\Gamma \times \varphi} A_j$, and $b:\prod_{\Gamma \times \psi} A$ with $\_ : (j^*b = a)$:
	\[\begin{tikzcd}
	& {\widetilde{\UU}} \\
	{\Gamma \times \psi} & \UU & {\Gamma \times \varphi}
	\arrow["\pi"', two heads, from=1-2, to=2-2]
	\arrow["b"{description}, from=2-1, to=1-2]
	\arrow["A"', from=2-1, to=2-2]
	\arrow[""{name=0, anchor=center, inner sep=0}, "a"{description}, curve={height=-6pt}, from=2-3, to=1-2]
	\arrow[""{name=1, anchor=center, inner sep=0}, "{j^*b}"{description}, shift right, curve={height=6pt}, from=2-3, to=1-2]
	\arrow["{j^*A}", from=2-3, to=2-2]
	\arrow[shift right=2, shorten <=4pt, shorten >=4pt, Rightarrow, no head, from=1, to=0]
\end{tikzcd}\]
\end{enumerate}
\end{proposition}

\begin{proof}
	For~(\ref{it:ext-form}), we use~Proposition~\ref{prop:dist-law} (and exponential transposition), finding:
	\[ (\Gamma \to \UU^{\ExtAt{j}}) = \left(\Gamma \to \sum_{A:\psi\to \UU} \prod_\varphi j^*A\right) = \sum_{A:\psi \times \Gamma \to \UU} \prod_{\Gamma \times \varphi} (\Gamma \times j)^* A\] 
	We use this and another application of~Proposition~\ref{prop:dist-law} to derive (\ref{it:ext-intro}):
	\begin{align*}
		(\Gamma \to \UU^{\LambdaAt{j}}) & \stackrel{\text{def.}}{=} \left( \Gamma \to \widetilde{\UU}^\psi \right) \\
		& \stackrel{\text{def.}}{=}  \left( \Gamma \to \sum_{A: \UU^\psi}  \prod_\psi A \right)
		\\ & = \left( \Gamma \to \sum_{A: \UU^\psi} \sum_{a:\prod_\varphi j^*A} \sum_{b: \prod_\psi A} (j^*b = a) \right) \\
		& = \sum_{A : \psi \times \Gamma \to \UU} \sum_{a : \prod_{\varphi \times \Gamma} (j \times \Gamma)^*A} \sum_{b : \prod_{\Gamma \times \psi} A} ( j^*b = a)
	\end{align*}
	Using the standard interpretation, it follows that these types correspond to the commutative diagrams as stated.
\end{proof}

\begin{proposition}
	The postulated bijections are \emph{natural} in the underlying context $\Gamma$, \ie, for any morphism $\sigma: \Delta \to \Gamma$ we have that the composite
\[\begin{tikzcd}
	\Delta && \Gamma && {\UU^\ExtAt{j}}
	\arrow["\sigma", from=1-1, to=1-3]
	\arrow["\alpha = \pair{A}{a}", from=1-3, to=1-5]
\end{tikzcd}\]
classifies pairs $\pair{\sigma^*A}{\sigma^*a}$ arising as in the diagram:
\[\begin{tikzcd}
	&& {\widetilde{\UU}} \\
	{\Delta \times \psi} & {\Delta \times \psi} & \UU
	\arrow["\pi", two heads, from=1-3, to=2-3]
	\arrow["{\sigma^*a}"{description}, from=2-1, to=1-3]
	\arrow["{\sigma \times \psi}"', from=2-1, to=2-2]
	\arrow["{\sigma^*A}"{description}, curve={height=24pt}, from=2-1, to=2-3]
	\arrow["a"{description}, from=2-2, to=1-3]
	\arrow["A"', from=2-2, to=2-3]
\end{tikzcd}\]
Given
\[\begin{tikzcd}
	\Delta && \Gamma && {\UU^\lambdaAt{j}}
	\arrow["\sigma", from=1-1, to=1-3]
	\arrow["\beta = \angled{A,a,b}", from=1-3, to=1-5]
\end{tikzcd}\]
yields in addition a section $b$ of $p_{A} : (\Gamma \times \psi).A \fibarr \Gamma \times \psi$, together with its reindexing $\sigma^*b$ such that all the following subdiagrams commute and the postulated equalities between the morphisms hold:
\[\begin{tikzcd}
	&& {\widetilde{\UU}} \\
	{\Delta \times \psi} & {\Gamma \times \psi} & \UU & {\Gamma \times \varphi} & {\Delta \times \varphi}
	\arrow["\pi", two heads, from=1-3, to=2-3]
	\arrow["{\sigma^*b}"{description}, from=2-1, to=1-3]
	\arrow["{\sigma \times \psi}"', from=2-1, to=2-2]
	\arrow["{\sigma^*A}"{description}, curve={height=24pt}, from=2-1, to=2-3]
	\arrow["b"{description}, from=2-2, to=1-3]
	\arrow["A"', from=2-2, to=2-3]
	\arrow[""{name=0, anchor=center, inner sep=0}, "a"{description, pos=0.4}, curve={height=-6pt}, from=2-4, to=1-3]
	\arrow[""{name=1, anchor=center, inner sep=0}, "{j^*b}"{description, pos=0.4}, curve={height=6pt}, from=2-4, to=1-3]
	\arrow["{j^*A}", from=2-4, to=2-3]
	\arrow[""{name=2, anchor=center, inner sep=0}, "{(j^*b)^\sigma}"{description, pos=0.4}, curve={height=28pt}, from=2-5, to=1-3]
	\arrow[""{name=3, anchor=center, inner sep=0}, "{\sigma^*a}"{description, pos=0.4}, curve={height=12pt}, from=2-5, to=1-3]
	\arrow["{(j^*A)^\sigma}"{description}, curve={height=-24pt}, from=2-5, to=2-3]
	\arrow["{\sigma \times \varphi}", from=2-5, to=2-4]
	\arrow[shift left, shorten <=2pt, shorten >=2pt, Rightarrow, no head, from=0, to=1]
	\arrow[shift left=2, shorten <=2pt, shorten >=2pt, Rightarrow, no head, from=3, to=2]
\end{tikzcd}\]
\end{proposition}

\begin{proof}
	This follows from~Propositions~\ref{prop:reindex-is-precomp} and \ref{prop:ext-class}.
\end{proof}

\theorem[Extension type structure]
Let $j \colon \varphi \hookrightarrow \psi$ be a shape inclusion.
\begin{enumerate}
	\item\label{it:ext-form-intro} The \emph{formation} and \emph{introduction rules} for extension types with respect to the shape inclusion $j$ are modeled by maps
	\begin{align*}
		\ExtAt{j} & \colon \UU^{\ExtAt{j}} \to \UU, \\
		\lambdaAt{j} & \colon \UU^{\lambdaAt{j}} \to \UU,
	\end{align*}
	making the following square commute:
\[\begin{tikzcd}
	{\UU^{\lambdaAt{j}}} && {\widetilde{\UU}} \\
	{\UU^{\ExtAt{j}}} && \UU
	\arrow["{\lambdaAt{j}}", from=1-1, to=1-3]
	\arrow[from=1-1, to=2-1]
	\arrow[from=1-3, to=2-3]
	\arrow["{\ExtAt{j}}"', from=2-1, to=2-3]
\end{tikzcd}\]
	\item\label{it:ext-elim-beta} Let us write $\UU^{\appAt{j}} := \UU^{\ExtAt{j}} \times_\UU \widetilde{\UU}$. The above square is a weak pullback, with a distinguished section $\appAt{j}$ as in
\[\begin{tikzcd}
	&& {\UU^{\lambdaAt{j}}} \\
	\\
	{\UU^{\appAt{j}} } && {\UU^{\appAt{j}} }
	\arrow[Rightarrow, no head, from=3-1, to=3-3]
	\arrow["\pi''", from=1-3, to=3-3]
	\arrow["{\appAt{j}}", from=3-1, to=1-3]
\end{tikzcd}\]
	if and only if the \emph{elimination rule with $\beta$-computation} holds.
	\item\label{it:ext-eta} It is a pullback if and only if in addition the \emph{$\eta$-computation rules} holds.
	\end{enumerate}

\endtheorem
	
\begin{proof}
	\textbf{Ad \ref{it:ext-form-intro}:} We have that
	\[ \UU^{\ExtAt{j}} = \sum_{A:\UU^\psi} \prod_\varphi j^* A \]
	and
	\[\UU^{\lambdaAt{j}} = \widetilde{\UU}^\psi = \sum_{A:\UU^\psi} \prod_\psi A =  \sum_{A:\UU^\psi}  \sum_{a : \Pi_\varphi j^*A} \sum_{b:\prod_\psi A} \big( j^*b = a\big).\]
	We have $\pi' : \UU^{\lambdaAt{j}} \to \UU^{\ExtAt{j}}$ with $\pi'\angled{A,a,b,*} := \angled{A,a}$. We directly see that the map
	\[ \ExtAt{j} : \UU^{\ExtAt{j}} \to \UU, \quad \angled{A,a} \mapsto \ExtAt{j}(A,a) \]
	models the introduction rule for extension types:
	\begin{mathpar}
		\inferrule*[right = (Ext-Form)]{\sh{t:I}{\varphi} \shape \\ \sh{t:I}{\psi} \shape \\ t:I \mid \varphi \types \psi \\
			\Xi\mid\Phi \types \Gamma \ctx \\
			\Xi,t:I \mid \Phi,\psi \mid \Gamma \types A\type \\
			\Xi,t:I \mid \Phi,\varphi \mid \Gamma \types a:A
		}{\Xi\mid\Phi\mid\Gamma \types \exten{t:I \mid \psi}{A}{\varphi}{a} \type}
	\end{mathpar}
	The map $\lambdaAt{j} : \UU^{\lambdaAt{j}} \to \widetilde{\UU}$ acts maps $\angled{A,a,b,*}$ into $\widetilde{\UU} = \sum_{T:\UU} \UU$. But commutativity of the diagram tells us that the type of the term that is the second component of $\lambdaAt{j}\angled{A,a,b,*}$ is $\ExtAt{j}(A,a)$. This means, $\lambdaAt{j}$ implements the introduction rule:
	\begin{mathpar}
		\inferrule*[right = (Ext-Intro)]{\sh{t:I}{\varphi} \shape \\ \sh{t:I}{\psi} \shape \\ t:I \mid \varphi \types \psi \\
			\Xi\mid\varphi \types \Gamma \ctx \\
			\Xi,t:I \mid \Phi,\psi \mid \Gamma \types A\type \\
			\Xi,t:I \mid \Phi,\varphi \mid \Gamma \types a:A \\\\
			\Xi,t:I \mid \Phi,\psi \mid \Gamma \types b:A \\
			\Xi,t:I \mid \Phi,\varphi \mid \Gamma \types b\jdeq a
		}{\Xi\mid\Phi\mid\Gamma \types \lam{t^{I\mid\psi}} b : \exten{t:I \mid \psi}{A}{\varphi}{a}}\and
	\end{mathpar}
	Given $A,a,b$ as in the premise, the $\lambda$-term $\lambda t^{I|\psi}.b : \exten{t:I|\psi}{A}{\varphi}{a}$ together with its typing is formed by composition:
	\[\begin{tikzcd}
\Gamma && {\UU^{\LambdaAt{j}}} && {\widetilde{\UU}}
\arrow["{\langle A,a,b\rangle}", from=1-1, to=1-3]
\arrow["{\lambdaAt{j}(A,a,b)}"{description}, curve={height=24pt}, from=1-1, to=1-5]
\arrow["{\lambdaAt{j}}", from=1-3, to=1-5]
\end{tikzcd}\]
\textbf{Ad \ref{it:ext-elim-beta}:}
We turn to interpreting the elimination rule  
	\begin{mathpar}	
		\inferrule*[right = (Ext-Elim)]{\sh{t:I}{\varphi} \shape \\ \sh{t:I}{\psi} \shape \\ t:I \mid \varphi \types \psi \\\\
			\Xi\mid\Phi \types \Gamma \ctx \\
			\Xi,t:I \mid \Phi,\psi \mid \Gamma \types A\type \\
			\Xi,t:I \mid \Phi,\varphi \mid \Gamma \types a:A
			\\\\
			\Xi\mid\Phi\mid\Gamma \types f:\exten{t:I \mid \psi}{A}{\varphi}{a} \\
			\Xi\types s:I \\ \Xi\mid\varphi\types \psi[s/t]
		}{\Xi\mid\Phi\mid\Gamma \types f(s) : A(s)}
	\end{mathpar}
	together with the computation rule
	\begin{mathpar}
		\inferrule*[right = (Ext-Comp)]{\sh{t:I}{\varphi} \shape \\ \sh{t:I}{\psi} \shape \\ t:I \mid \varphi \types \psi \\\\
	\Xi\mid\Phi \types \Gamma \ctx \\
	\Xi,t:I \mid \Phi,\psi \mid \Gamma \types A\type \\
	\Xi,t:I \mid \Phi,\varphi \mid \Gamma \types a:A
	\\\\
	\Xi\mid\Phi\mid\Gamma \types f:\exten{t:I \mid \psi}{A}{\varphi}{a} \\
	\Xi\types s:I \\ \Xi\mid\Phi\types \psi[s/t]
}{\Xi\mid\Phi\mid\Gamma \types f(s) \jdeq a[s/t] : A(s)}
\end{mathpar}
	and the $\beta$-rule:
	\begin{mathpar}
		\inferrule*[right = (Ext-$\beta$)]{\sh{t:I}{\varphi} \shape \\ \sh{t:I}{\psi} \shape \\ t:I \mid \varphi \types \psi \\
	\Xi\mid\Phi \types \Gamma \ctx \\
	\Xi,t:I \mid \Phi,\psi \mid \Gamma \types A\type \\
	\Xi,t:I \mid \Phi,\varphi \mid \Gamma \types a:A \\\\
	\Xi,t:I \mid \Phi,\psi \mid \Gamma \types b:A \\
	\Xi,t:I \mid \Phi,\varphi \mid \Gamma \types b\jdeq a\\
	\Xi\types s:I \\ \Xi\mid\Phi\types \psi[s/t]
}{\Xi\mid\Phi\mid\Gamma \types (\lam{t^{I\mid\psi}} b)(s) \jdeq b[s/t]}
	\end{mathpar}
	For a family $A : \Gamma \times \Psi \to \UU$ and a section $f : \Gamma \times \Psi \to \widetilde{\UU}$ with $\pi \circ f = A$ we model $A[s]$ and $b[s]$ by precomposition with $s : 1 \to \Gamma \times \Psi$:
\[\begin{tikzcd}
&& {\widetilde{\UU}} \\
1 & {\Gamma \times \Psi} \\
&& \UU
\arrow["s", from=2-1, to=2-2]
\arrow["b", from=2-2, to=1-3]
\arrow["A"', from=2-2, to=3-3]
\arrow["\pi", from=1-3, to=3-3]
\arrow["{b[s]}"{description}, curve={height=-12pt}, from=2-1, to=1-3]
\arrow["{A[s]}"{description}, curve={height=12pt}, from=2-1, to=3-3]
\end{tikzcd}\]
We assume $\UU^{\lambdaAt{j}}$ is a weak pullback of $\pi : \widetilde{\UU} \to \UU$ and $\ExtAt{j} : \UU^{\ExtAt{j}} \to \UU$. Then the premise of the elimination rule gives us:
\[\begin{tikzcd}
\Gamma \\
& {\UU^{\LambdaAt{j}}} && {\widetilde{\UU}} \\
& {\UU^{\ExtAt{j}}} && \UU
\arrow[from=2-2, to=2-4]
\arrow[from=2-2, to=3-2]
\arrow["\pi", from=2-4, to=3-4]
\arrow["{\ExtAt{j}}"', from=3-2, to=3-4]
\arrow["f", curve={height=-12pt}, from=1-1, to=2-4]
\arrow["{\langle A,a\rangle}"', curve={height=12pt}, from=1-1, to=3-2]
\arrow["{\langle A,a, \widetilde{f}\rangle}"{pos=0.9}, dashed, from=1-1, to=2-2]
\end{tikzcd}\]
By our description of $\UU^{\LambdaAt{j}}$ a map $\langle A,a, \widetilde{f}\rangle : \Gamma \to \UU^{\LambdaAt{j}}$ corresponds to a section $\psi \; | \; \Gamma \vdash f : A$ that strictly restricts to $a$ over $\varphi$:
\[\begin{tikzcd}
& {\widetilde{\UU}} \\
{\Gamma \times \psi} & \UU
\arrow["A"', from=2-1, to=2-2]
\arrow["\pi", from=1-2, to=2-2]
\arrow["{\widetilde{f}}", from=2-1, to=1-2]
\end{tikzcd}\]
\[\begin{tikzcd}
{\Gamma \times \varphi} & {\Gamma \times \psi} & {\widetilde{\UU}}
\arrow["{\widetilde{f}}"', from=1-2, to=1-3]
\arrow["a", curve={height=-18pt}, from=1-1, to=1-3]
\arrow["{\Gamma \times j}"', from=1-1, to=1-2]
\end{tikzcd}\]
We then set
\[ f(s) := \widetilde{f}[s] = \widetilde{f} \circ s,\]
so that $\psi \; | \; \Gamma \vdash f(s) : A[s]$ since:
\[\begin{tikzcd}
&& {\widetilde{\UU}} \\
1 & {\Gamma \times \Psi} & \UU
\arrow["A"', from=2-2, to=2-3]
\arrow["\pi", from=1-3, to=2-3]
\arrow["{{\widetilde{f}}}", from=2-2, to=1-3]
\arrow["{f(s)}", curve={height=-18pt}, from=2-1, to=1-3]
\arrow["s", from=2-1, to=2-2]
\arrow["{A[s]}"', curve={height=24pt}, from=2-1, to=2-3]
\end{tikzcd}\]
For the $\beta$-computation rule, let us consider a map $\angled{A,a,b} : \Gamma \to \UU^{\LambdaAt{j}}$. Consider the induced map in the following pullback:
\[\begin{tikzcd}
\Gamma \\
& {\UU^{\appAt{j}}} && {\widetilde{\UU}} \\
& {\UU^{\ExtAt{j}}} && \UU
\arrow["{{\langle A,a,\widetilde{\lambdaAt{j}(A,a,b)}\rangle}}"{description}, dashed, from=1-1, to=2-2]
\arrow["{\lambdaAt{j}_{A,a}(b)}", curve={height=-12pt}, from=1-1, to=2-4]
\arrow["{{\pair{A}{a}}}"', curve={height=34pt}, from=1-1, to=3-2]
\arrow["{\lambdaAt{j}}", from=2-2, to=2-4]
\arrow[two heads, from=2-2, to=3-2]
\arrow["\lrcorner"{anchor=center, pos=0.125}, draw=none, from=2-2, to=3-4]
\arrow["\pi", two heads, from=2-4, to=3-4]
\arrow["{{\ExtAt{j}}}"', from=3-2, to=3-4]
\end{tikzcd}\]
But coming back to the initial data
\[\begin{tikzcd}
&& {\widetilde{\UU}} \\
{\Gamma \times \varphi} & {\Gamma \times \psi} & \UU
\arrow["\pi", two heads, from=1-3, to=2-3]
\arrow["a", curve={height=-6pt}, from=2-1, to=1-3]
\arrow["{\Gamma \times j}"', tail, from=2-1, to=2-2]
\arrow["b", from=2-2, to=1-3]
\arrow["A"', from=2-2, to=2-3]
\end{tikzcd}\]
we see that this induces the map $\angled{A,a,b} : \Gamma \to \UU^{\appAt{j}}$ that also is a gap map for the above pullback, hence $\Gamma | \psi \vdash b \jdeq \widetilde{\lambdaAt{j}(A,a,b)}$. Thus, for $s : \psi$ we get
\[ \Gamma \vdash (\lambdaAt{j}_{A,a}(b))(s) = (\widetilde{\lambdaAt{j}_{A,a}(b)})[s] = b[s] : A[s].\]

\textbf{Ad \ref{it:ext-eta}:} For the $\eta$-computation rule, assume that $\UU^{\LambdaAt{j}}$ is an actual, strict pullback of $\widetilde{\UU}$ and $\UU^{\ExtAt{j}}$.

We first establish some more setup. Given a family $A : \Gamma \times \psi.A \to \UU$, consider its reindexing along $p_A : \Gamma \times \psi.A \to A$, namely the family $Ap_A : \Gamma\times \psi.A \to \UU$ given as follows:
\[\begin{tikzcd}
	{\Gamma\times \psi.A.Ap_{A}} & {\Gamma\times \psi.A} & {\widetilde{\UU}} \\
	{\Gamma \times \psi.A} & {\Gamma \times \psi} & \UU
	\arrow["{q_{Ap_A}}", two heads, from=1-1, to=1-2]
	\arrow["{p_{Ap_A}}"', two heads, from=1-1, to=2-1]
	\arrow["\lrcorner"{anchor=center, pos=0.125}, draw=none, from=1-1, to=2-2]
	\arrow["{q_A}", from=1-2, to=1-3]
	\arrow["{p_A}"', two heads, from=1-2, to=2-2]
	\arrow["\lrcorner"{anchor=center, pos=0.125}, draw=none, from=1-2, to=2-3]
	\arrow["\pi", two heads, from=1-3, to=2-3]
	\arrow["{p_A}", two heads, from=2-1, to=2-2]
	\arrow["{A{p_A}}"', curve={height=18pt}, from=2-1, to=2-3]
	\arrow["A", from=2-2, to=2-3]
\end{tikzcd}\]
We will abbreviate $A' := Ap_A : \Gamma \times \psi.A \to \UU$.

Then from
\[\begin{tikzcd}
	{(\Gamma \times \psi).A} & {\widetilde{\UU}} \\
	{\Gamma \times \psi} & \UU
	\arrow["{q_A}", from=1-1, to=1-2]
	\arrow["{p_A}"', two heads, from=1-1, to=2-1]
	\arrow["\lrcorner"{anchor=center, pos=0.125}, draw=none, from=1-1, to=2-2]
	\arrow["\pi", two heads, from=1-2, to=2-2]
	\arrow["A"', from=2-1, to=2-2]
\end{tikzcd}\]
we get that $q_A$ can be viewed as a section of the family $A'$:
\[\begin{tikzcd}
	{\Gamma \times \psi.A.Ap_A} & {\widetilde{\UU}} \\
	{\Gamma \times \psi.A} & \UU
	\arrow["{{q_{A'}}}", from=1-1, to=1-2]
	\arrow["{{p_{A'}}}", two heads, from=1-1, to=2-1]
	\arrow["\lrcorner"{anchor=center, pos=0.125}, draw=none, from=1-1, to=2-2]
	\arrow["\pi", two heads, from=1-2, to=2-2]
	\arrow["{{q_A}}", curve={height=-18pt}, dotted, from=2-1, to=1-1]
	\arrow["{{A'}}", from=2-1, to=2-2]
\end{tikzcd}\]
Furthermore, let $a : \Gamma \times \varphi \to j^*A$ be a partial section of $A$. We can then consider its pullback along $p_{j^*A}$, written $a' := p_{j^*A}^*(a)$:
\[\begin{tikzcd}
	{\Gamma \times \varphi.j^*A.j^*Ap_{j^*A}} && {\Gamma \times \varphi.j^*A} \\
	{\Gamma \times \varphi.j^*A} && {\Gamma \times \varphi}
	\arrow[from=1-1, to=1-3]
	\arrow["{q_{j^*Ap_{j^*A}}}", from=1-1, to=2-1]
	\arrow["\lrcorner"{anchor=center, pos=0.125}, draw=none, from=1-1, to=2-3]
	\arrow["{p_{j^*A}}", from=1-3, to=2-3]
	\arrow["{a'}", curve={height=-18pt}, dotted, from=2-1, to=1-1]
	\arrow["{p_{j^*A}}", from=2-1, to=2-3]
	\arrow["a", curve={height=-18pt}, dotted, from=2-3, to=1-3]
\end{tikzcd}\]
Let us abbreviate $A'' := j^*Ap_{j^*A} : \Gamma \times \varphi.j^*A \to \UU$.

Now, consider a term $f$ of $\ExtAt{j}(A,a)$:
\[\begin{tikzcd}
	\Gamma \\
	& {\UU^{\LambdaAt{j}}} && {\widetilde{\UU}} \\
	& {\UU^{\ExtAt{j}}} && \UU
	\arrow["{{\langle A,a,\widetilde{f}\rangle}}"{pos=0.8}, dashed, from=1-1, to=2-2]
	\arrow["f", curve={height=-12pt}, from=1-1, to=2-4]
	\arrow["{{\langle A,a \rangle}}"', curve={height=12pt}, from=1-1, to=3-2]
	\arrow["{{\LambdaAt{j}}}", from=2-2, to=2-4]
	\arrow[two heads, from=2-2, to=3-2]
	\arrow["\pi", two heads, from=2-4, to=3-4]
	\arrow["{\ExtAt{j}}", from=3-2, to=3-4]
\end{tikzcd}\]
Now, precomposition by $p_A$ yields the term
\[\begin{tikzcd}
	{\Gamma \times \psi.A} & \Gamma & {\widetilde{\UU}}
	\arrow["{p_A}", from=1-1, to=1-2]
	\arrow["{fp_A}"{description}, curve={height=18pt}, from=1-1, to=1-3]
	\arrow["f", from=1-2, to=1-3]
\end{tikzcd}\]
fitting into the diagram:
\[\begin{tikzcd}
	{\Gamma \times \psi.A} \\
	& {\UU^{\LambdaAt{j}}} && {\widetilde{\UU}} \\
	& {\UU^{\ExtAt{j}}} && \UU
	\arrow["{\langle A',a',\widetilde{fp_A}\rangle}"{pos=0.8}, dashed, from=1-1, to=2-2]
	\arrow["{fp_A}", curve={height=-12pt}, from=1-1, to=2-4]
	\arrow["{\langle A',a' \rangle}"', curve={height=12pt}, from=1-1, to=3-2]
	\arrow["{\LambdaAt{j}}", from=2-2, to=2-4]
	\arrow[two heads, from=2-2, to=3-2]
	\arrow["\lrcorner"{anchor=center, pos=0.125}, draw=none, from=2-2, to=3-4]
	\arrow["\pi", two heads, from=2-4, to=3-4]
	\arrow["{\ExtAt{j}}"', from=3-2, to=3-4]
\end{tikzcd}\]
This allows us to consider ``application'' of $fp_A$ to $q_A$, giving rise to the term $fp_A(q_A) : \Gamma \times \psi \to \widetilde{\UU}$:
\[\begin{tikzcd}
	&&&& {\widetilde{\UU}} \\
	{\Gamma \times \psi} && {\Gamma \times \psi.A.Ap_A} && \UU
	\arrow["\pi", two heads, from=1-5, to=2-5]
	\arrow["{(fp_A)(q_A)}", curve={height=-18pt}, from=2-1, to=1-5]
	\arrow["{\langle1,q_A\rangle}", from=2-1, to=2-3]
	\arrow["{\widetilde{fp_A}}", from=2-3, to=1-5]
	\arrow["{A'}", from=2-3, to=2-5]
\end{tikzcd}\]
We note that $\widetilde{fp_A} = \widetilde{f}p_A$. Thus
\[ (fp_A)(q_A) = \widetilde{fp_A}[q_A] = \widetilde{fp_A}(1,p_A) = \widetilde{f}p_A(1,q_A)=\widetilde{f}.\]
From this, we get
\[ \lambdaAt{j}_{A,a}(fp_A(q_A)) = \lambdaAt{j} \circ (A,a,fp_A(q_A)) = \lambdaAt{j} \circ \widetilde{f} = f, \]
as desired.
\end{proof}

The strict stability result can then be phrased as the following theorem:
\theorem[Strict stability of extension types]\label{thm:coh-ext}
	Let $\E$ be a type-theoretic model topos, and consider the interpretation of simplicial homotopy type theory \`{a} la \cite{RS17,Shu19,Riehl-Sem} in the type-theoretic model topos $\sE$. Then, a fixed splitting of the standard type formers induced by the universal fibration $\pi: \widetilde{\UU} \to \UU$ for small fibrations yields an interpretation of the extension types $\ExtAt{j}$ (where $j:\varphi \to_I \psi$ is a fixed shape inclusion relative to some cube $I$) that is strictly stable under pullback. This means, for families $A: \Gamma \times \psi \to \UU$ and partial sections $a:\Gamma \times \varphi \to_{\Gamma \times \psi} A$ we have: 
	\begin{enumerate}
		\item\label{it:stab-i} $\sigma^*\big(\ExtAt{j}(A,a)\big) = \ExtAt{j}(\sigma^*A,\sigma^*a)$
		\item\label{it:stab-ii} $\sigma^*\big(\LambdaAt{j}_{A,a}(b)\big) = \LambdaAt{j}_{\sigma^*A,\sigma^*a}(\sigma^*b)$
		\item\label{it:stab-iii} $\sigma^*\big(\appAt{j}_{A,a}(f,s)\big) = \appAt{j}_{\sigma^*A,\sigma^*a}(\sigma^*f,\sigma^*s)$
	\end{enumerate}
\endtheorem
In what follows, we establish the proof of this theorem by splitting the extension type formers, assuming all the necessary and pre-established logical and model-categorical structure in the background.

\begin{proof}
	\textbf{Eq.~\ref{it:stab-i}:} The formation and introduction rule are, top to bottom, Rule 1 and 2, respectively, in~\cite[Figure~4]{RS17} (or \cf~Figure \ref{fig:exten-i}).

	Morphisms from $\Gamma$ into the object $\UU^{\LambdaAt{j}}$ consist of triples $\angled{A,a,b}$ as in
	\[\begin{tikzcd}
		&& A \\
		\\
		{\Gamma \times \varphi} && {\Gamma \times \psi}
		\arrow["{\Gamma \times j}"', from=3-1, to=3-3]
		\arrow[two heads, from=1-3, to=3-3]
		\arrow["a", dotted, from=3-1, to=1-3]
		\arrow["b"', curve={height=18pt}, dotted, from=3-3, to=1-3]
	\end{tikzcd}\]
	\ie, strictly $a = b \circ (\Gamma \times j)$. Considering the projection $\UU^{\LambdaAt{j}} \to \UU^{\ExtAt{j}}$ the fiber at an instance $\pair{A}{a}$ is exactly the semantic extension type
	\[ \ExtAt{j}(A,a):=\ccexten{\psi}{A}{\varphi}{a} \]
	defined by the \emph{split} cartesian square, after~\cite[Theorem~A.16]{RS17}. This means, the object $\ccexten{\psi}{A}{\varphi}{a}$ and its projection to $\Gamma$ have been \emph{chosen}:
	\[\begin{tikzcd}
		\ccexten{\psi}{A}{\varphi}{a}  && {A^\psi} \\
		\\
		\Gamma && {A^\varphi \times_{(\Gamma \times \psi)^\varphi} (\Gamma \times \psi)^\psi}
		\arrow["{\pair{a}{\eta} }", from=3-1, to=3-3]
		\arrow["\lrcorner"{anchor=center, pos=0.125}, draw=none, from=1-1, to=3-3]
		\arrow[two heads, from=1-1, to=3-1]
		\arrow[two heads, from=1-3, to=3-3]
		\arrow[from=1-1, to=1-3]
	\end{tikzcd}\]
	In particular, the right vertical map is a fibration, so the object $\ExtAt{j}(A,a)$ is fibrant. If shapes are taken to be fibrant, said map being a fibration follows from type-theoretic-ness of the model structure. Otherwise one would have to use that the model structure is cartesian monoidal as in~\cite[Lemma~A.4]{RS17}.
	
	Note that in the extensional type theory of the ambient topos one could describe the extension type also as the $\Sigma$-type
	\[ \sum_{b:\psi \to A} (b\circ j \jdeq a),\]
	where $(-_1) \jdeq (-_2)$ stands for the extensional identity type.
	
	Now, in analogy to \cite[Theorem~1.4.15]{KV18}, we define the map $\ExtAt{j}:\UU^{\lambdaAt{j}} \to \widetilde{\UU}$ by the universal property of the extension type. The generic context $\UU^{\lambdaAt{j}}$ consists of triples $\angled{A,a,b}$~s.t.
	\[\begin{tikzcd}
		&& {A^\psi} \\
		\\
		\Gamma && {A^\varphi \times_{{\Gamma \times \psi}^\varphi} (\Gamma \times \psi)^\psi}
		\arrow["b", from=3-1, to=1-3]
		\arrow["{\pair{a}{\eta} }"', from=3-1, to=3-3]
		\arrow[two heads, from=1-3, to=3-3]
	\end{tikzcd}\]
	where $\eta: \Gamma \to (\Gamma \times \psi)^\psi$ denotes the transpose of the identity of $\Gamma \times \psi$. Now, we define the generic $\lambda$-term of the extension type as the gap map $(\LambdaAt{j})(b) := \widehat{b}$ of the pullback:
	\[\begin{tikzcd}
		\Gamma \\
		&& \ccexten{\psi}{A}{\varphi}{a} && {A^\psi} \\
		\\
		&& \Gamma && {A^\varphi \times_{{\Gamma \times \psi}^\varphi} (\Gamma \times \psi)^\psi}
		\arrow[two heads, from=2-5, to=4-5]
		\arrow[from=2-3, to=2-5]
		\arrow[two heads, from=2-3, to=4-3]
		\arrow["{\pair{a}{\eta} }", from=4-3, to=4-5]
		\arrow["b", curve={height=-18pt}, from=1-1, to=2-5]
		\arrow["{\widehat{b}}"{description}, dashed, from=1-1, to=2-3]
		\arrow[curve={height=18pt}, Rightarrow, no head, from=1-1, to=4-3]
		\arrow["\lrcorner"{anchor=center, pos=0.125}, draw=none, from=2-3, to=4-5]
	\end{tikzcd}\]
	As shown prior, defining
	\[\begin{tikzcd}
		{\UU^{\lambdaAt{j}} := \widetilde{\UU}^\psi} \\
		&& {\UU^{\ExtAt{j}} := \UU^\psi \times_{\UU^\varphi} \widetilde{\UU}^\varphi} &&&& {\widetilde{\UU}^\varphi} \\
		\\
		&& {\UU^\psi} &&&& {\UU^\varphi}
		\arrow["{\UU^j}"', from=4-3, to=4-7]
		\arrow["{\pi^\varphi}", from=2-7, to=4-7]
		\arrow[from=2-3, to=4-3]
		\arrow[from=2-3, to=2-7]
		\arrow["\lrcorner"{anchor=center, pos=0.125}, draw=none, from=2-3, to=4-7]
		\arrow["{\pi^\psi}"', curve={height=24pt}, from=1-1, to=4-3]
		\arrow["{\widetilde{\UU}^j}", curve={height=-24pt}, from=1-1, to=2-7]
		\arrow[from=1-1, to=2-3]
	\end{tikzcd}\]
	then the square
	\[\begin{tikzcd}
		{\UU^\LambdaAt{j}} &&& \widetilde{\UU} \\
		\\
		{\UU^\ExtAt{j}} &&& \UU
		\arrow[from=1-1, to=3-1]
		\arrow["{\ExtAt{j}}"', from=3-1, to=3-4]
		\arrow["{\LambdaAt{j}}", from=1-1, to=1-4]
		\arrow[from=1-4, to=3-4]
		\arrow["\lrcorner"{anchor=center, pos=0.125}, draw=none, from=1-1, to=3-4]
		\arrow["\lrcorner"{anchor=center, pos=0.125}, draw=none, from=1-1, to=3-4]
	\end{tikzcd}\]
	being a pullback precisely captures the universal property.
	
	In the terminology of \cite[Theorem~2.3.4]{KL18} this means that the projection $\UU^{\LambdaAt{j}} \to \UU^{\ExtAt{j}}$ is the \emph{universal dependent extension type} over $\UU^{\ExtAt{j}}$. In particular, in presence of the splitting it is a \emph{chosen} small fibration. Thus, as illustrated \eg~in \cite[Remark~2.6]{Awo18} we obtain strict pullback stability: everything in
	\[\begin{tikzcd}
		&& \Delta \\
		\\
		\Gamma && {\UU^{\ExtAt{j}}} && \UU
		\arrow["{\pair{A}{a}}"{description}, from=3-1, to=3-3]
		\arrow["{\ExtAt{j}}"{description}, from=3-3, to=3-5]
		\arrow["\sigma"', from=1-3, to=3-1]
		\arrow["{(\ExtAt{j})(\sigma^*A,\sigma^*a)}", from=1-3, to=3-5]
		\arrow["{\pair{\sigma^*A}{\sigma^*a}}"{description}, from=1-3, to=3-3]
		\arrow["{(\ExtAt{j})(A,a)}"{description}, curve={height=24pt}, from=3-1, to=3-5]
	\end{tikzcd}\]
	commutes strictly yielding as desired
	\[ \ExtAt{j}(A,a) \circ \sigma = \ExtAt{j}(A\sigma, a\sigma). \]
	\textbf{Eq.~\ref{it:stab-ii}:}

	Similarly, and using that the generic lifts $b$ are given by gap maps of (strict) pullbacks we obtain
	\[ \lambdaAt{j}(A,a,b) \circ \sigma = \ExtAt{j}(A\sigma, a\sigma, b\sigma). \]
	\textbf{Eq.~\ref{it:stab-iii}:}

	The elimination and computation rule are Rule 3 and 4, respectively, in~\cite[Figure~4]{RS17} (or \cf~Figure~\ref{fig:exten-ii}).

	To interpet the elimination rule, note first that the exponential $A^\psi$ comes with a (chosen) evaluation map:
	\[\begin{tikzcd}
		{A^\psi \times_\Gamma (\Gamma \times \psi)} && A \\
		& \Gamma
		\arrow["{\mathrm{ev}}", from=1-1, to=1-3]
		\arrow[from=1-1, to=2-2]
		\arrow[from=1-3, to=2-2]
	\end{tikzcd}\]
	By pulling this back along the (chosen) map from the extension type, we obtain a map
	\[ \mathrm{ev}: \ExtAt{j}(A,a) \times_\Gamma (\Gamma \times \psi) \to_\Gamma A,\]
	fibered over $\Gamma$. Now,we want to define application $\appAt{j}$ of a function $f:\Gamma \to_\Gamma \ExtAt{j}(A,a)$ to a section $s: \Gamma \to_\Gamma \Gamma \times \psi$. Analogously to~\cite[Theorem~1.4.15]{KL18}, we take the composition:
	\[\begin{tikzcd}
		\Gamma && {\ExtAt{j}(A,a) \times_\Gamma (\Gamma \times \psi)} && A \\
		&& \Gamma
		\arrow[from=1-3, to=2-3]
		\arrow["{\langle f,s \rangle}", from=1-1, to=1-3]
		\arrow["{\mathrm{ev}}", from=1-3, to=1-5]
		\arrow[Rightarrow, no head, from=1-1, to=2-3]
		\arrow[from=1-5, to=2-3]
		\arrow["{\appAt{j}(f,s)}", curve={height=-40pt}, from=1-1, to=1-5]
	\end{tikzcd}\]
	This validates the rule, since $\appAt{j}(f,s):\Gamma \to_\Gamma A$, as desired.
	For substitution along $\sigma: \Delta \to \Gamma$, consider the following diagram involving chosen fibrations and split cartesian squares:
	\[\begin{tikzcd}
		{\ExtAt{j} (\sigma^*A,\sigma^*a) \mathrlap{\times_\Delta(\Delta \times \psi)}} &&& {\ExtAt{j} (A,a) \times_\Gamma(\Gamma \times \psi)} \\
		& {\sigma^*A} &&& A \\
		\Delta &&& \Gamma & {}
		\arrow[from=1-1, to=3-1]
		\arrow["{\appAt{j}(\sigma^*f,\sigma^*s)}"{pos=0.65}, shorten <=4pt, dashed, from=1-1, to=2-2,]
		\arrow[shorten <=56pt, from=1-1, to=1-4]
		\arrow[from=2-2, to=3-1]
		\arrow["{\pair{\sigma^*f}{\sigma^*a}}", curve={height=-18pt}, dotted, from=3-1, to=1-1]
		\arrow[from=3-1, to=3-4, ]
		\arrow[from=1-4, to=3-4]
		\arrow["{\appAt{j}(f,s)}"{description}, from=1-4, to=2-5]
		\arrow[from=3-4, to=2-5]
		\arrow["\lrcorner"{anchor=center, pos=0.125}, shift left=2, draw=none, from=1-1, to=3-4]
		\arrow["\lrcorner"{anchor=center, pos=0.125}, draw=none, from=2-2, to=3-4]
		\arrow["{\pair{f}{a}}"{description, pos=0.7}, curve={height=-18pt}, dotted, from=3-4, to=1-4]
		\arrow[from=2-2, to=2-5, crossing over]
	\end{tikzcd}\]
	This uniquely defines the map $\appAt{j}(\sigma^*f,\sigma^*s)$, and by construction
	\[ \sigma^*\appAt{j}(f,s) = \appAt{j}(\sigma^*f,\sigma^*s). \]
	Furthermore the desired computation rule holds, saying that for a term $s$ in the smaller tope $\varphi$, one judgmentally has $\appAt{j}(f,s) \jdeq a[s]$. This is established by the commutation of the diagram:
	\[\begin{tikzcd}
		&& A \\
		\\
		{\Gamma \times \varphi} && {\Gamma \times \psi} \\
		& \Gamma
		\arrow["{\Gamma \times j}"{description}, from=3-1, to=3-3]
		\arrow[from=1-3, to=3-3]
		\arrow["a"{description}, from=3-1, to=1-3]
		\arrow["s"{description}, from=4-2, to=3-1]
		\arrow[from=4-2, to=3-3]
		\arrow["f"{description}, curve={height=-12pt}, dotted, from=3-3, to=1-3]
		\arrow["{\appAt{j}(f,s)}"', curve={height=70pt}, dotted, from=4-2, to=1-3]
	\end{tikzcd}\]
	
	The $\beta$- and $\eta$-rule (\cf~\cite[Figure~4]{RS17}) can be proven to hold similarly (Rule 5 and 6, respectively, in~\cite[Figure~4]{RS17}, or \cf~Figure~\ref{fig:exten-iii}).
\end{proof}

\section{Outlook}

We have used Veovdosky's splitting method from~\cite[Section~4, Definition~0.1]{VVTySys} to split the extension type formers from Riehl--Shulman's synthetic $\inftyone$-category theory~\cite{RS17}. Our presentation also closely follows and parallels Awodey's work on natural models~\cite{Awo18}.

Notably, we have only worked using a single universe. We believe there to be no essential problems generalizing this line of reasoning to the case of a hierarchy of universes, \cf~\cite{LW15,Shu19,RS17}. We also worked in a simplified technical setting with fibrant shapes, blurring the lines to the pre-type layers which are treated quite explicitly in~\cite[Appendix~A]{RS17}. While this can be justified semantically, one could still try to adapt the splitting to a setting more faithful to the original theory~\cite{RS17} and its fibered structure/multi-part contexts. Already in~\cite{RS17} it is suggested that it should be possible to do so along the lines of~\cite{LW15}.

Yet another line of investigation could consist of trying to capture extension types in the recent $2$-categorical framework of Coraglia--Di~Liberti~\cite{CDL-Ctx}.

\appendix

\section{Rules for extension types}\label{app:exten}

The rules for the extension types of sHoTT have been taken from~\cite[Figure~4]{RS17}. In simplicial HoTT, a general typing judgment has the form
\[ \Xi\;|\;\varphi\;|\;\Gamma \vdash A \; \mathrm{type}\]
where $\Xi$ denotes a cube context, $\varphi$ a shape context, and $\Gamma$ a type context in the standard sense. The rules are analogous to that of $\Pi$-types, but with strict judgmental side conditions added.

\begin{figure}[h]
	\centering
	\begin{mathpar}
		\inferrule*[right = (Ext-Form)]{\sh{t:I}{\varphi} \shape \\ \sh{t:I}{\psi} \shape \\ t:I \mid \varphi \types \psi \\
			\Xi\mid\varphi \types \Gamma \ctx \\
			\Xi,t:I \mid \Phi,\psi \mid \Gamma \types A\type \\
			\Xi,t:I \mid \Phi,\varphi \mid \Gamma \types a:A
		}{\Xi\mid\Phi\mid\Gamma \types \exten{t:I \mid \psi}{A}{\varphi}{a} \type}\and
		\inferrule*[right = (Ext-Intro)]{\sh{t:I}{\varphi} \shape \\ \sh{t:I}{\psi} \shape \\ t:I \mid \varphi \types \psi \\
			\Xi\mid\Phi \types \Gamma \ctx \\
			\Xi,t:I \mid \Phi,\psi \mid \Gamma \types A\type \\
			\Xi,t:I \mid \Phi,\varphi \mid \Gamma \types a:A \\\\
			\Xi,t:I \mid \Phi,\psi \mid \Gamma \types b:A \\
			\Xi,t:I \mid \Phi,\varphi \mid \Gamma \types b\jdeq a
		}{\Xi\mid\Phi\mid\Gamma \types \lam{t^{I\mid\psi}} b : \exten{t:I \mid \psi}{A}{\varphi}{a}}
	\end{mathpar}
	\caption{Extension types: formation and introduction}
	\label{fig:exten-i}
\end{figure}

\begin{figure}[h]
	\centering
	\begin{mathpar}	
		\inferrule*[right = (Ext-Elim)]{\sh{t:I}{\varphi} \shape \\ \sh{t:I}{\psi} \shape \\ t:I \mid \varphi \types \psi \\\\
			\Xi\mid\Phi \types \Gamma \ctx \\
			\Xi,t:I \mid \Phi,\psi \mid \Gamma \types A\type \\
			\Xi,t:I \mid \Phi,\varphi \mid \Gamma \types a:A
			\\\\
			\Xi\mid\Phi\mid\Gamma \types f:\exten{t:I \mid \psi}{A}{\varphi}{a} \\
			\Xi\types s:I \\ \Xi\mid\Phi\types \psi[s/t]
		}{\Xi\mid\Phi\mid\Gamma \types f(s) : A(s)}\and
		\inferrule*[right = (Ext-Comp)]{\sh{t:I}{\varphi} \shape \\ \sh{t:I}{\psi} \shape \\ t:I \mid \varphi \types \psi \\\\
			\Xi\mid\Phi \types \Gamma \ctx \\
			\Xi,t:I \mid \Phi,\psi \mid \Gamma \types A\type \\
			\Xi,t:I \mid \Phi,\varphi \mid \Gamma \types a:A
			\\\\
			\Xi\mid\Phi\mid\Gamma \types f:\exten{t:I \mid \psi}{A}{\varphi}{a} \\
			\Xi\types s:I \\ \Xi\mid\varphi\types \varphi[s/t]
		}{\Xi\mid\Phi\mid\Gamma \types f(s) \jdeq a[s/t] : A(s)}\and
	\end{mathpar}
	\caption{Extension types: elimination and computation}
	\label{fig:exten-ii}
\end{figure}

\begin{figure}[h]
		\centering
		\begin{mathpar}	
		\inferrule*[right = (Ext-$\beta$)]{\sh{t:I}{\varphi} \shape \\ \sh{t:I}{\psi} \shape \\ t:I \mid \varphi \types \psi \\
			\Xi\mid\Phi \types \Gamma \ctx \\
			\Xi,t:I \mid \Phi,\psi \mid \Gamma \types A\type \\
			\Xi,t:I \mid \Phi,\varphi \mid \Gamma \types a:A \\\\
			\Xi,t:I \mid \Phi,\psi \mid \Gamma \types b:A \\
			\Xi,t:I \mid \Phi,\varphi \mid \Gamma \types b\jdeq a\\
			\Xi\types s:I \\ \Xi\mid\varphi\types \psi[s/t]
		}{\Xi\mid\Phi\mid\Gamma \types (\lam{t^{I\mid\psi}} b)(s) \jdeq b[s/t]}\and  
		\inferrule*[right = (Ext-$\eta$)]{\sh{t:I}{\varphi} \shape \\ \sh{t:I}{\psi} \shape \\ t:I \mid \varphi \types \psi \\\\
			\Xi\mid \Phi \types \Gamma \ctx \\
			\Xi,t:I \mid \Phi,\psi \mid \Gamma \types A\type \\
			\Xi,t:I \mid \Phi,\varphi \mid \Gamma \types a:A
		\\\\
			\Xi\mid\Phi\mid\Gamma \types f:\exten{t:I \mid \psi}{A}{\varphi}{a} \\
		}{\Xi\mid\Phi\mid\Gamma \types f \jdeq \lam{t^{I\mid\psi}} f(t)}\and
	\end{mathpar}
	\caption{Extension types: $\beta$-conversion and $\eta$-abstraction}
	\label{fig:exten-iii}
\end{figure}

\end{document}